\documentclass[12pt,reqno]{amsart}
\usepackage{amsmath,amssymb}
\usepackage[dvipdfmx]{graphicx,xcolor}
\usepackage{tikz} 
\usetikzlibrary{positioning, intersections, calc, arrows.meta,math} 
\usepackage{latexsym,amssymb}
\usepackage{mathrsfs}
\usepackage[abbrev]{amsrefs}
\usepackage{comment}

\newtheorem{thm}{Theorem}[section]
\newtheorem{cor}[thm]{Corollaly}

\newtheorem{lem}[thm]{Lemma}

\renewcommand{\th}{\theta}

\newcommand{\ep}{\varepsilon}

\renewcommand{\div}{{\rm {div\ }}}

\def\<{\langle }
\renewcommand{\>}{\rangle }

\newcommand{\supp}{\text{ supp }}

\newcommand{\bl}{\color{blue}}
\renewcommand{\bl}{\color{black}}

 \newenvironment{pf}
    {{\noindent \bf Proof. }}{\hfill $\Box$}

\numberwithin{equation}{section}
\numberwithin{thm}{section}

\begin{document}

\begin{center}\large \bf 
Energy conservation law for weak solutions
of the full compressible Navier-Stokes equations
\end{center}

\footnote[0]
{

\indent
E-mail: 
\, $^*$motofumi.aoki.q2@dc.tohoku.ac.jp
$^{**}$t-iwabuchi@tohoku.ac.jp
}
\vskip5mm

\begin{center}
{\sf 
Motofumi Aoki$^{*}$
Tsukasa Iwabuchi$^{**}$
}

\vskip2mm

Mathematical Institute, 
Tohoku University\\
Sendai 980-8578 Japan

\end{center}

\vskip5mm

\begin{center}
\begin{minipage}{135mm}
\footnotesize
{\sc Abstract. }
We consider a sufficient condition for the energy conservation law of a weak solution for the full compressible Navier-Stokes equations on the torus. 
We prove that a weak solution constructed by Feireisl 
with certain integrability conditions conserve the energy. 
Our assumption relaxes the regularity condition 
for the density compared with existing results.  
\end{minipage}
\end{center}

\section{Introduction}

\subsection{The Equations}\mbox{}

The full system of compressible Navier-Stokes equations in $\mathbb T^d $ with $ d \ge 2$ is written as follows. 
\begin{equation}
\label{cNS1}
\left\{
\begin{split}
&\partial_t \rho + {\rm div \,}(\rho u) =0,
\hspace{5mm} &t > 0, \hspace{1.2mm} x \in \mathbb T^d , \\
	&\partial_t (\rho u) + {\rm div}\,( \rho u \otimes u)  + \nabla p( \rho, \theta) - {\rm div}\, \mathbb S = 0,
	\hspace{5mm} &t > 0, \hspace{1.2mm} x \in \mathbb T^d , \\
	&\partial_t \left(\rho \left(\frac{|u|^2}{2} + 
	e(\rho, \theta)\right)\right) + 
	{\rm div \,}\left(\rho u \left(\frac{|u|^2}{2} + e(\rho, \theta) \right) \right) 
	\\&\hspace{36mm}+{\rm div\,} q = 
	{\rm div\,}(\mathbb S u) - {\rm div\,}(p(\rho, \theta) u)
	\hspace{5mm} &t > 0, \hspace{1.2mm} x \in \mathbb T^d , \\
	&( \rho, u, e(\rho, \theta)) \big|_{t=0} = (\rho_0, u_0, e(\rho_0, \theta_0)),\hspace{5mm} &x \in \mathbb T^d,  
\end{split}
\right.
\end{equation}
where $\rho = \rho(t,x)$, $u = (u_1(t,x), \cdots , u_d (t,x))$, $p$ and $e(\rho, \theta)$ denote the unknown density of the fluid, the unknown velocity vector, the unknown pressure of the fluid and the unknown internal energy of the fluid at the point $(t,x) \in (0,T) \times \mathbb T^d$.  
$\rho_0 = \rho_0 (x)$ is the initial density, 
$u_0 = (u_{0,1}(x), \cdots, u_{0,d}(x))$ is the initial velocity vector 
and $e(\rho_0, \theta_0) = e(\rho_0(x), \theta_0 (x))$ is the initial internal energy. 
$\mathbb S, \mu, \lambda$ satisfies 
\[
\begin{cases}
\mathbb S := \mu \{ \nabla u + (\nabla u)^T\} + \lambda \, {\rm div}\, u \, \mathbb I,  \\                                                                                                                         \mu > 0, \,\,
\displaystyle
\lambda +\frac{2}{d}\mu \ge 0.
\end{cases}
\]
The equations \eqref{cNS1} consist of the continuity, the motion, and the total energy of the fluid. 
The purpose of this paper is to study a sufficient condition to obtain 
the convervation of the total energy.

Let us start by a simple argument to obtain the energy conservation for the equations \eqref{cNS1}. 
We focus on weak solutions such that 
\begin{equation}\label{0124-1}
\rho \in L^{\infty}(0,T; L^{\infty}), 
\quad u \in L^{\infty}(0,T; L^2) \cap L^2(0,T; W^{1,2}), 
\quad e(\rho, \theta) \in L^{\infty}(0,T; L^1). 
\end{equation}
If we consider the equations in the sense of the distributions, then 
the following integrability condition,  
\begin{equation} \notag 
\rho u^2, e (\rho, \theta), \rho u |u|^2, \rho e(\rho, \theta) u , q, \mathbb S \, u , p(\rho, \theta) 
 \in L^1(0,T; L^1),
\end{equation}
is sufficient to have the energy conservation that 
\[
\int_{\mathbb T^d} \rho(t) \Big( \frac{|u(t)|^2}{2} + e(\rho(t) , \theta(t))\Big) dx 
= \int_{\mathbb T^d} \rho_0 \Big( \frac{|u_0|^2}{2} + e(\rho_0 , \theta_0)\Big) dx .
\]
By this argument, we naturally have the following proposition.

\vskip3mm 


\noindent
{\bf Proposition.} Let $d \ge 2$ and
$( \rho , u, e(\rho, \theta))$ be a weak solution of \eqref{cNS1} 
which satisfies \eqref{0124-1} and 
\[
\begin{aligned}
u \in L^3(0,T ;L^3), \ 
e(\rho, \theta) \in L^{\infty}(0,T; L^2), \ 
q, \, p (\rho, \theta) \in L^1(0,T ;L^1). 
\end{aligned}
\]
Then the solution of \eqref{cNS1} satisfies the energy equality.

\vskip3mm 

\noindent
{\bf Remark.} In the case when $d = 2, 3, 4$, 
$u \in L^\infty (0,T ; L ^2) \cap L^2 (0, T; W^{1,2})$ implies the 
integrability $u \in L^3 (0,T ; L^3)$. 

\vskip3mm 

\subsection{Existing Results}\mbox{}

To the best of our knowledge, 
there is not any result of the energy conservation for the full-system of compressible Navier-Stokes 
equations~\eqref{cNS1}, and our motivation comes from 
existing papers. 
In this subsection, 
we refer to results of the energy conservation of weak solutions 
for the incompressible Navier-Stokes equations, 
inhomogeneous incompressible Navier-Stokes equations, 
and compressible Navier-Stokes equations. 
We also mention several papers related to Onsagar's conjecture 
to understand the difference between 
the viscous and the inviscid cases. 

\vskip3mm 
 
\noindent
{\bf A. Navier-Stokes equations.}

\begin{enumerate}
\item (The incompressbile Navier-Stokes) \, 
Leray \cite{Le-34} first showed the energy identity for a weak solution 
in the two dimensional case. 
In the higher dimensional case $d \ge 3$, 
we need to suppose some regularity condition to obtain the energy identity. 
Serrin~\cite{Se-62} showed the regularity of solutions in $L^p(0,T;L^q)$ 
with $2/p + d/q \le 1$ and $q \ge d$, which leads to the energy identity. 
Later Shinbrot~\cite{Sh-74} studied in a larger class 
$L^p(0,T;L^q)$ with $ 2/p + 2/q \le 1$ and $q \ge 4$. 
The most general results are those in 
$L^{p,w}(0,T; B^{\frac{2}{p}+\frac{2}{q}-1}_{q, \infty})$ with $1 \le p < q \le \infty$, $1/p + 2/q <1$ 
and in 
$L^p(0,T; B^{\frac{5}{2p}+\frac{3}{q}-\frac{3}{2}}_{q,\infty})$ 
with $0 < p \le 3$, $1 \le q \le \infty$, $1/p + 2/q \le 1$, 
by Cheskidov--Luo~\cite{ChLu-20}.

\item (Inhomogeneous incompressbile Navier-Stokes) \, 
Leslie--Shvydkoy \cite{LeSh-16} first showed the energy identity 
for a weak solution such that 
$u \in L^p(0,T;B^{\frac{1}{3}}_{p,c_0}(\mathbb T^3)),\, \rho \in L^q(0,T; B^{\frac{1}{3}}_{q, \infty} (\mathbb T^3))$, $p \in L^{\frac{p}{2}}(0,T; B^{\frac{1}{3}}_{\frac{b}{2}, \infty})$, $3/p + 1/q =1$, 
and with some assumption for the pressure. 
Recently Nguyen--Nguyen--Tang \cite{NgNgTa-19} considered the same framework 
for the velocity as in Shinbrot~\cite{Sh-74}.

\item (Barotropic compressible Navier-Stokes) \, 
Yu \cite{Yu-17} proved the energy identity for isentropic Navier-Stokes equations 
if the velocity belongs to $L^p(0,T; L^q(\mathbb T^3))$ with $ 2/p + 2/q \le 5/6$ and $q \ge 6$. 
Akramov--D\polhk ebiec--Skipper--Wiedemann \cite{AkDeSkWi-20} showed the energy identity when   
$u \in B^{\alpha}_{3, \infty}((0,T) \times \mathbb T^3) \cap L^2(0,T;W^{1,2})$, $\rho, \rho u \in B^{\beta}_{3, \infty}((0,T) \times \mathbb T^3)$ with $\alpha + 2\beta >1,\,  2 \alpha + \beta >1,\, 0 \le \alpha, \beta \le 1$ and 
the pressure is continuous with respect to the density. 
Furthermore, Nguyen--Nguyen--Tang
\cite{NgNgTa-19} 
considered weak solutions such that 
$u \in L^p(0,T;L^q), \, \rho \in L^{\infty}(0,T; B^{\frac{1}{2}}_{\frac{12}{5},\infty}(\mathbb T^3))$ with $2/p + 2/q=1, \, q \ge 4$, 
and the pressure is a $C^2$ function with respect to the density,  
to have the energy identity. 
\end{enumerate}

As described above, we have the energy identity of weak solutions if we impose 
some regularity condition. 
In general, the energy identity of weak solutions have been left as an open problem.


\vskip3mm 

\noindent
{\bf B. Onsager Conjecture's View.}

Onsager \cite{On-49} conjectured that the energy of solutions 
for the incompressible Euler equations, 
in the H\"older spaces with exponent $\alpha$, 
is conserved when $\alpha > 1/3$ 
and that the energy is possibly dissipated when $\alpha < 1/3$.
Eyink \cite{Ey-94} and Constantin--E.--Titi~\cite{CoETi-94} first proved that weak solutions $u \in L^3(0,T; B^{\alpha}_{3, \infty})$ with $\alpha >1/3$ conserve the energy. 
Cheskidov--Constantin--Friedlander--Shvydkoy \cite{ChCoFrSh-08}
generalized the condition to $\alpha = 1/3$. 
Analogous result for the incompressible Navier-Stokes equations 
was studied by Cheskidov--Luo~\cite{ChLu-20}. 


We mention several results for the other equations. 
For the inhomogeneous Euler equations, Feireisl--Gwiazda--Gwiazda--Wiedemann \cite{FeGwGwWi-17} proved the energy identity local in time for solutions such that 
$u \in B^{\alpha}_{p, \infty}((0,T) \times \mathbb T^3)$, $\rho, \rho u \in B^{\beta}_{q, \infty}((0,T) \times \mathbb T^3)$, $p \in L^{p'}_{loc}((0,T) \times \mathbb T^3)$
 with $1 \le p,q \le \infty, 0 \le \alpha, \beta \le 1$, $2/p + 1/q =1, \, 1/p + 1/p' =1$ and $2\alpha + \beta >1$. 
They also studied barotropic Euler equations where $\rho, u, \rho u$ are in the same class, $p \in C^2$ and $\alpha +2\beta > 1$.
This assumption is generalized by Akramov--D\polhk ebiec--Skipper--Wiedemann \cite{AkDeSkWi-20} to the pressure $p \in C^{1,(\gamma -1)}$, $\alpha +\gamma \beta >1$ and $1 < \gamma < 2$.
Drivas--Eyink~\cite{DrEy-18} proved the energy conservation law for full Euler system such that $u \in B^{\alpha, \infty}_{p,loc}((0,T) \times \mathbb T^3), \rho \in B^{\beta, \infty}_{p,loc} ((0,T) \times \mathbb T^3), \theta \in B^{\gamma, \infty}_{p,loc}((0,T) \times \mathbb T^3)$ with 
$2\min\{\beta, \gamma\} + \alpha >1, \min\{\beta, \gamma\} + 2\alpha >1, 3 \alpha >1$ and $p \ge 3$

\vskip3mm

\noindent {\bf C. Strong solutions. }

There are many studies of the Cauchy problem for \eqref{cNS1} under the scaling critical setting.
Danchin~\cite{Da-011}, \cite{Da-012} proved the well-posedness with the initial data  $(\rho_0, u_0, \theta_0) \in (\dot B^{\frac{d}{p}}_{p,1}, \dot B^{-1+\frac{d}{p}}_{p,1}, \dot B^{-2+\frac{d}{p}}_{p,1})$ with $1 < p < d$.
Furthermore Chikami--Danchin \cite{ChDa-15} discuss the unique solvability with $(\rho_0 -1, u_0, E_0) \in (\dot B^{\frac{d}{p}}_{p,1}, \dot B^{-1+\frac{d}{p}}_{p,1}, \dot B^{-2+\frac{d}{p}}_{p,1})$ where $1 < p < 2d$ and  $E_0 := \displaystyle \frac{\rho_0 u^2_0}{2} + e(\rho_0, \theta_0)$.

On the other hand, there are some studies about the ill-posedness for the motion of the ideal gas. 
Chen--Miao--Zhang \cite{ChMiZh-15} proved the ill-posedness in the case when $d=3$ for the initial data $(\rho_0, u_0, \theta_0) \in (\dot B^{\frac{3}{p}}_{p,1}, \dot B^{-1+\frac{3}{p}}_{p,1}, \dot B^{-2+\frac{3}{p}}_{p,1}) $ with $p>3$.
Recently, Iwabuchi and Ogawa \cite{IwOg-21} proved the ill-posedness in the case when $d=2$ for the initial data $(\rho_0, u_0, \theta_0) \in (\dot B^{\frac{2}{p}}_{p,q}, \dot B^{-1+\frac{2}{p}}_{p,q}, \dot B^{-2+\frac{2}{p}}_{p,q})$ with $1 \le p \le \infty$ and $1 \le q < \infty$.

\subsection{Main Results}\mbox{}

By assuming a relation for the elastic potential and the pressure,  
we rewrite the equations~\eqref{cNS1} following Feireisl~\cite{Fe-04}. 
The difference of the momentum equation multiplied by $u$ 
and the continuity equation multiplied by $|u|^2/2$ is 
the following equation for the kinetic energy. 
\[
\partial_{t} \left(\rho \frac{|u|^2}{2} \right) + 
{\rm div\,} \left(\rho u \frac{|u|^2}{2} \right) = 
({\rm div \,} \mathbb S) \cdot u + 
(\nabla p(\rho, \theta)) \cdot u.
\]
We notice that the total energy is the sum of the kinetic energy 
and the internal energy $e (\rho, \theta)$, 
which consists of the elastic potential $h(\rho)$
and the thermal energy contribution $Q(\theta)$. 
Then we can write the following equation for the internal energy. 
\begin{equation}
\label{0201-1}
\partial_t (\rho e(\rho, \theta)) + 
{\rm div\,} (\rho e(\rho, \theta) u) + {\rm div \,} q = 
(\mathbb S : \nabla u) + p(\rho, \theta)\, {\rm div\,} u.
\end{equation}
Here $(\mathbb S : \nabla u)$ denotes a product of matrices, i.e. for matrix $A =\{a_{ij}\}, B = \{b_{ij}\}$, 
\[
\displaystyle
(A: B) := \bigg\{ \sum_{1\le i,j \le d} \{a_{ij} b_{ji} \}_{ij} \bigg\}.
\] 
\noindent
Next we decompose \eqref{0201-1} into the two energies due to the elastic potential $h(\rho)$ and due to the thermal energy contribution $Q(\theta)$. 
We suppose the internal energy $e(\rho, \theta)$ and the pressure 
$p(\rho, \theta)$, 
\[
\begin{cases}
e(\rho, \theta) := h(\rho) + Q(\theta), \\
p(\rho, \theta) := p_e(\rho) + p_{th} (\rho, \theta) = p_e(\rho) + \theta \partial_{\theta} p(\rho),
\end{cases}
\]
where the elastic potential and the thermal energy contribution 
are given by 
\begin{equation}\label{0202-2}
\begin{cases}
\displaystyle
h(\rho) := \int^{\rho}_1 \frac{p_e(z)}{z^2} dz, \\
\displaystyle
Q( \theta) := \int^{\theta}_0 c_v(z) dz, \quad c_v(z) \ge c_v > 0, \, z \ge 0.
\end{cases}
\end{equation}
By multiplying the equation of continuity by $(\rho h(\rho))'$, we get
\begin{equation}
\label{1026-3}
\partial_t(\rho h(\rho)) + {\rm div \,} (\rho h(\rho) u) + 
p_e(\rho) {\rm div\,} u=0. 
\end{equation}
We subtract \eqref{1026-3} from \eqref{0201-1}, 
and then obtain the equation of the thermal energy $Q(\theta)$.  
\[
\partial_t (\rho Q(\theta)) + {\rm div}\, (\rho Q(\theta) u) - {\rm div \,} q = (\mathbb S : \nabla u) - p_{th} ( \rho, \theta) {\rm div \,} u.
\]
Moreover, $q$ is given by Fourier's law. 
\[
q := - \nabla \mathcal K(\theta). 
\]
Thus we are led to the following equations. 
\begin{equation} \label{cNS3}
\left\{
	\begin{split}
	&\partial_t \rho + {\rm div \,} (\rho u) =0, 
	&t > 0, \, x \in \mathbb T^d , \\
	&\partial_t (\rho u)  + {\rm div}\,( \rho u \otimes u)  + \nabla p( \rho, \theta) - {\rm div}\, \mathbb S = 0,
	&t > 0, \, x \in \mathbb T^d , \\
	&\partial_t (\rho Q(\theta)) + {\rm div}\, (\rho Q(\theta) u) - \Delta \mathcal K (\theta) = (\mathbb S : \nabla u) - p_{th} ( \rho, \theta) {\rm div \,} u,
	&t > 0, \, x \in \mathbb T^d , \\
&( \rho, u, \theta)\big|_{t=0} = (\rho_0, u_0, \theta_0),
&x \in \mathbb T^d , 
	\end{split}
\right.	
	\end{equation}

The aim of this paper is to give a sufficient condition that weak solutions satisfy the energy equality, 
\[
E(t) = E(0), \quad t > 0,
\]
where 
\begin{equation}\label{0202-5}
\displaystyle
E(t) := 
\frac{1}{2} \int_{\mathbb T^d} \rho u^2 (t,x) dx  + 
\int_{\mathbb T^d} \rho h (\rho) (t,x) dx + \int_{\mathbb T^d}\rho Q(\theta) (t,x) dx. 
\end{equation}
We introduce a weak solution of \eqref{cNS3} based on an idea of Feireisl~\cite{Fe-04}. 

 
\noindent 
{\bf Definition}.
$\Omega \subset \mathbb R^d$, $d \ge 2$, is bounded domain. 
Let $T > 0$ and $\displaystyle \gamma > \frac{d}{2}$. 
A measurable function $(\rho, u, \theta)$ on $(0,T) \times \Omega $ is called a weak solution of \eqref{cNS3} on $(0,T)$ if 
\begin{enumerate}
\item $\rho, u, \theta$ satisfy for $\displaystyle r > \frac{2d}{d+2}$,\\
$
\rho \in L^{\infty}(0,T; L^{\gamma}), \, 
u \in L^2 (0,T; W^{1,2}_0), \, 
\rho u \in L^{\infty}(0,T; L^{\frac{2 \gamma}{\gamma +1}}), \, 
\mathbb S  \in L^2(0,T; L^2), \\ 
p_e(\rho) \in L^1(0,T;L^1), \,
\rho Q(\theta) \in L^{\infty}(0,T; L^1) \cap L^2(0,T ; L^2 \cap L^r), \, \mathcal K(\theta) \in L^1(0,T; L^1),\\
\log \, \theta \in L^2(0,T; L^2).$
\item $\rho,\, u$ satisfy the equation of continuty in the distribution sense, i.e.
\[
\begin{aligned}
\int^t_0 \{ \< \rho, \partial_{\tau} \Phi \> \} + \< \rho u, \nabla \Phi \> \} d\tau = \<\rho(t), \Phi(t)\>-\<\rho_0, \Phi_0\>, 
\end{aligned}
\]
for any test function $\Phi \in \mathcal D([0,T) \times \Omega)$ with $\Phi > 0$.
\item $\rho,\, u$ satisfy the momentum equation in the distribution sense, i.e.
\[
\begin{aligned}
\int^t_0 \{ \< \rho u , \partial_{\tau} \phi \> + 
\< \rho u \otimes u, \nabla \phi \> &+ 
\< p(\rho, \theta) , {\rm div}\, \phi \> \} d\tau -
\int^t_0 \< \mathbb S, \nabla \phi\> d\tau  \\
&=
\< \rho u(t), \phi(t)\> -
\< \rho_0 u_0, \phi (0)\> ,
 \end{aligned}
\]
for any test function $\phi \in \mathcal D([0,T) \times \Omega)^d$.
\item $\rho,\, u, \, \theta$ satisfy the internal equation in the distribution sense, i.e.
\[
\begin{aligned}
 \int^t_0 \{ &\< \rho Q(\theta), \partial_{\tau} \psi \> + 
 \< \rho Q(\theta) u , \nabla \psi \> - 
 \< Q(\theta) , \Delta \psi \> \} d\tau \\
&=
\<\rho(t)Q(\theta)(t), \psi(t) \> - 
\< \rho_0Q(\theta_0), \psi(0) \> -
\int^t_0 \{\< (\mathbb S , \nabla u), \psi \> - \< \theta \partial_{\theta} p \, {\rm div \,} u, \psi \> \}d\tau, \,\,
\end{aligned}
\]
\noindent
for any test function $\psi \in \mathcal D([0,T) \times \Omega)$.
\item $\rho, \,\rho u$ and $\rho Q(\theta)$ satisfy the initial conditions
\[
\begin{cases}
\< \rho(t, x) , \eta (x)\> &\to \< \rho_0, \eta(x) \>,\\
 \< \rho (t,x) u (t,x) , \eta(x) \> &\to  \< \rho_0 u_0 , \eta(x) \>, \\
 \< \rho(t, x) Q(\theta)(t,x), \eta (x)\> &\to
 \< \rho_0 Q(\theta_0), \eta (x)\>,
\end{cases}
\]
for any tset function $\eta \in \mathcal D(\Omega)$.
\end{enumerate}

\vskip3mm 


This is our main result of this paper. 



\begin{thm}
\label{1129-2}
Let $d =2,3$, $\Omega = \mathbb T^d$. Suppose that $T <\infty$, $(\rho, u, \theta)$ 
is a weak solution of \eqref{cNS3} on $(0,T)$.
Assume that 
\begin{equation}
\begin{aligned}
\label{0915-3}
&0 < c_1 \le \rho \le c_2 < \infty, \, \quad 
u \in L^{\infty}(0,T; L^2), 
\end{aligned}
\end{equation}
for some $c_1, c_2$, and $p \in C^1([c_1 , c_2])$.  
In the case when $d=3$, 
we additionally assume that 
\begin{equation}
\label{0915-5}
\begin{split}
\displaystyle 
& 
u \in L^p(0,T; L^q), 
\quad \frac{1}{p} + \frac{3}{q} = 1,
\quad 3 \le q \le 4,
\\
\displaystyle
\text{or } \quad 
& u \in L^p(0,T; L^q),
\quad \frac{2}{p} + \frac{2}{q}= 1, 
\quad q \ge 4. 
\end{split}
\end{equation}
Then 
the energy \eqref{0202-5} conserves, i.e., 
$E(t) = E(0)$ for all $0 < t < T$. 
\end{thm}

It is also possible to apply the condition \eqref{0915-5} to the barotropic compressible Navier-Stokes equations.
\begin{equation} \label{cNS4}
\left\{
	\begin{split}
	&\partial_t \rho + {\rm div \,} (\rho u) =0, 
	&t > 0, \, x \in \mathbb T^d , \\
	&\partial_t (\rho u)  + {\rm div}\,( \rho u \otimes u)  + \nabla p( \rho) - {\rm div}\, \mathbb S = 0,
	&t > 0, \, x \in \mathbb T^d , \\
&( \rho, u)\big|_{t=0} = (\rho_0, u_0),
&x \in \mathbb T^d. 
	\end{split}
\right.	
	\end{equation}

\begin{cor}
Let $d =2,3$, $\Omega = \mathbb T^d$. Suppose that $T <\infty$, $(\rho, u)$ 
is a weak solution of \eqref{cNS4} on $(0,T)$, and 
\[
\begin{aligned}
0 < c_1 \le \rho \le c_2 < \infty, \quad 
u \in L^{\infty}(0,T; L^2),   
\end{aligned}
\]
for some $c_1, c_2$, and $p \in C^1([c_1 ,c_2])$.
In the case when $d=3$, 
we additionally assume \eqref{0915-5}. 
Then the energy equality holds. 
\[
\begin{aligned}
& \displaystyle
\frac{1}{2} \int_{\mathbb T^d}  \rho  u^2 (t,x) dx +
\int_{\mathbb T^d} \rho h(\rho) (t,x) dx  
\\
& + 
\mu \int^t_0 \int_{\mathbb T^d} | \nabla u + (\nabla u)^{T}|^2 dx d\tau + 
\lambda \int^t_0 \int_{\mathbb T^d} |{\rm div \,} u|^2 dx d\tau \\
=&  
\frac{1}{2} \int_{\mathbb T^d}  \rho_0  u^2_0 (x) dx +
\int_{\mathbb T^d} \rho_0 h(\rho_0) (x) dx 
\end{aligned}
\]
for all $0< t < T$, where 
\[
h (\rho) : = 
\int^{\rho}_1 \frac{p(z)}{z^2} dz.
\]

\end{cor}

\noindent
{\bf Remark}. 
\begin{enumerate}
\item In the case when $d=2$, we have by the Gagliardo-Nirenbarg inequality 
and the Ho\"older inequality that 
\[
\| u\|_{L^p(0,T;L^q)}
\le 
C\|u\|^{1-\frac{2}{p}}_{L^{\infty}(0,T; L^2)}
\|u\|^{\frac{2}{p}}_{L^2(0,T;W^{1,2})}, 
\]
for $2/p + 2/q =1$ with $2\le q < \infty$.
We do not need the additional assumption~\eqref{0915-5}. 

\item  In the case when $d=3$, 
the condition $u \in L^4(0,T;L^4)$ is the best in \eqref{0915-5}, 
and the other cases follow from the following inequalities 
shown by the H\"older interpolation and the emebedding $W^{1,2} \subset L^6$. 
\[
\begin{aligned}
 \, &\| u\|_{L^4(0,T;L^4)}
\le \| u \|^{1- \frac{p}{2(p-2)}}_{L^2(0,T;W^{1,2})}
\|u\|^{\frac{p}{2(p-2)}}_{L^p(0,T;L^q)}
, 
\quad 
&\frac{1}{p}+ \frac{3}{q} =1, \quad \dfrac{1}{4} \le \dfrac{1}{q} &\le \dfrac{1}{3},
\\
 \, &\| u\|_{L^4(0,T;L^4)}
\le \| u \|^{1- \frac{p}{4}}_{L^{\infty}(0,T;L^2)}
\|u\|^{\frac{p}{4}}_{L^p(0,T;L^q)}, 
\quad 
&\frac{2}{p}+ \frac{2}{q} =1, \,\,\,\,\, \quad \quad \dfrac{1}{q} &\leq \dfrac{1}{4}.
\end{aligned}
\]


\item Comparing with the previous studies on incompressible fluids, our result corresponds to the result by Shinbrot \cite{Sh-74}.

\item 
We have obtained the energy identity without the positive regularity assumption 
for the density and the pressure by Nguyen--Nguyen--Tang \cite{NgNgTa-19}, 
where they suppose 
$\rho \in L^{\infty}(0,T; B^{\frac{1}{2}}_{2, \infty})$ when $d=2$ and $\rho \in L^{\infty}(0,T;B^{\frac{1}{2}}_{\frac{12}{5}, \infty})$ when $d=3$. 
Moreover we have generalized the assumption $p \in C^2(0,\infty)$ to $p \in C^1(0,\infty)$.

\item Akramov--D\polhk ebiec--Skipper--Wiedemann \cite{AkDeSkWi-20} showed the energy identity for the compressible Navier-Stokes equations such that $u \in B^{\alpha}_{3, \infty}((0,T) \times \mathbb T^3) \cap L^2(0,T;W^{1,2})$, $\rho, \rho u \in B^{\beta}_{3, \infty}((0,T) \times \mathbb T^3)$, $0 \le c_1 \le \rho \le c_2 < \infty$ with $\alpha + 2\beta >1, 2 \alpha + \beta >1, 0 \le \alpha, \beta \le 1$ and $p \in C([c_1, c_2])$. 
Our theorem does not require the positive regularity assumption in the Besov space, 
while we assume  $p \in C^1 ([c_1, c_2])$.

\end{enumerate}

\vskip3mm 

 We finally mention the energy equality for a compressible Navier-Stokes equations which describe the motion of the ideal gas.
Taking the pressure as $p(\rho, \theta)= \rho \theta$ and 
the internal energy $Q(\theta) = \theta$ in \eqref{cNS3}, 
we write the motion of the ideal gas as follows.
\begin{equation} \label{cNS5}
\left\{
	\begin{split}
	&\partial_t \rho + {\rm div \,} (\rho u) =0, 
	&t > 0, \, x \in \mathbb T^d , \\
	&\partial_t (\rho u)  + {\rm div}\,( \rho u \otimes u)  + \nabla (\rho \theta) - {\rm div}\, \mathbb S = 0,
	&t > 0, \, x \in \mathbb T^d , \\
	&\partial_t (\rho \theta) + {\rm div}\, (\rho \theta u) - \Delta \theta = (\mathbb S : \nabla u) - \rho \theta {\rm div \,} u,
	&t > 0, \, x \in \mathbb T^d , \\
&( \rho, u, \theta)\big|_{t=0} = (\rho_0, u_0, \theta_0),
&x \in \mathbb T^d. 
	\end{split}
\right.	
	\end{equation}

We can also lead to the regularity condition which \eqref{cNS5} satisfies the energy identity from Theorem \ref{1129-2};

\begin{cor}
\label{1130-3}
Let $d =2,3$ and $\Omega = \mathbb T^d$. Suppose that $T <\infty$, $(\rho, u, \theta)$ 
is a weak solution of \eqref{cNS5} on $(0,T)$.
Assume that 
\[
\begin{aligned}
\label{1130-1}
&0 < c_1 \le \rho \le c_2 < \infty, \, 
u \in L^{\infty}(0,T; L^2).
\end{aligned}
\]
for some $c_1, c_2$. 
In the case when $d=3$ we additionally assume \eqref{0915-5}. 
Then the energy equality holds, i.e., 
\[
\displaystyle
\frac{1}{2} \int_{\mathbb T^d} \rho u^2 (t,x) dx  + 
 \int_{\mathbb T^d}\rho \theta (t,x) dx
 =
\frac{1}{2} \int_{\mathbb T^d} \rho_0 u_0^2 (t,x) dx  + 
 \int_{\mathbb T^d}\rho_0 \theta_0 (t,x) dx 
\]
for all $0< t < T$. 
\end{cor}

\noindent
{\bf Remark.} 
In the case when $d=2$, Iwabuchi and Ogawa \cite{IwOg-21} showed the norm inflation 
of smooth solutions to \eqref{cNS5} for the initial data in 
$\dot B^{\frac{2}{p}}_{p, q}\times \dot B^{-1+\frac{2}{p}}_{p,q}\times  \dot B^{-2 +\frac{2}{p}}_{p,q}$ with $1 \le p \le \infty$ and $1 \le q < \infty$.
On the other hand, a natural class for the energy is 
$L ^\infty \times L^2 \times L^1$. 
However, there is no inclusion relation between the two framework 
and it would be interesting to reveal the relation.

  

\vskip3mm

\noindent
{\bf Notations. }  
A test function space $\mathcal D (\mathbb T^d)$ is written by 
\[
\begin{aligned}
\mathcal D (\mathbb T^d) &: = 
\{ f \in C^{\infty}(\mathbb T^d) \big| \, \,
p_m(f) < \infty \,\,{\rm for \,\, all \,\,}
m \in \mathbb N \cup \{0\}
\},\\
p_m(f) &:= 
\sum_{|\alpha| \le m}
\sup_{x \in \mathbb T^d} | \partial^{\alpha} f(x) |, \quad
m = 0,1,2, \cdots.
\end{aligned}
\]
A test function space $\mathcal D=\mathcal D ([0,T) \times \mathbb T^d)$ is written by 
\begin{gather*}
\mathcal D : = 
\{ f \in C^{\infty}([0,T) \times \mathbb T^d) \big| \,
\supp f  \subset [0,T) \times \mathbb T^d, \,\,
p_m(f) < \infty \,\,{\rm for \,\, all \,\,}
m \in \mathbb N \cup \{0\}
\},\\
p_m(f) := 
\sum_{\alpha + |\beta| \le m}
\sup_{(t,x) \in [0,T) \times \mathbb T^d} | \partial^{\alpha}_t \partial^{\beta}_x f(x) |, \quad
m = 0,1,2, \cdots.
\end{gather*}


\section{Preliminary}


To prove Theorem \ref{1129-2}, we introduce a mollifier for space and a lemma.
\medskip 

\noindent
{\bf Definition}.
Let $\eta \in C^{\infty}_0(\mathbb R^{d})$ be such that 
$\eta$ is radially symmetric and 
\[
{\rm supp}\, \eta \subset B_1(0), \quad  0 \leq \eta \leq 1,
\quad \int_{\mathbb R^d} \eta(x) dx = 1.
\]
For $\varepsilon > 0$ we set 
\[
\eta_{\varepsilon}(x) = \varepsilon ^{-d} \eta(\varepsilon^{-1}x),  \quad 
x \in \mathbb T^d , 
\]
and define $u_{\varepsilon}(x)$ by
\[
u_{\varepsilon}(x)= (\eta_{\varepsilon} * u)(x) :=
 \int_{\mathbb T^d} \eta_{\varepsilon} (x-y)  u(y) dy, \quad 
 x \in \mathbb T^d .
\]

\begin{lem}{\rm (\cite{NgNgTa-19})}
\label{0909-1}
Let $ d\ge 2, 1 \le p,q \le \infty$. \\
\begin{enumerate}
%
\item 
There exists $C > 0$ such that for every $\ep > 0$ 
and every $f \in L^p(0,T; L^q(\mathbb T^d))$ 
\[
\| \nabla f_{\ep}\|_{L^p(0,T;L^{q})} \le C \ep^{-1} \|f\|_{L^p(0,T;L^q)}.
\]
Moreover, if $p,q <\infty$, then 
\begin{equation}\label{1208-1}
\limsup_{\ep \to 0} \ep \|\nabla f_{\ep}\|_{L^p(0,T;L^q)} = 0,
\end{equation}
provided that $f \in L^p(0,T; L^q(\mathbb T^d))$. 
\item Let $ c > 0$. Then 
for every $g \in L^\infty (0,T ; L^\infty)$ with $\inf g > c $,  
there exists $C > 0$ such that 
for every $f \in L^p(0,T; L^q) $ 
and every $\ep >0$
\[
\bigg\| \nabla \frac{f_{\ep}}{g_{\ep}}\bigg\|_{L^p(0,T;L^q(\mathbb T^d))} 
\le C \ep^{-1} \|f\|_{L^p(0,T;L^q(\mathbb T^d))}.
\]
Moreover, if $p,q <\infty$, then
\begin{equation}\label{1208-2}
\limsup_{\ep \to 0} \ep \bigg\| \nabla \frac{f_{\ep}}{g_{\ep}}\bigg\|_{L^p(0,T;L^q(\mathbb T^d))} = 0
\end{equation}
provided that 
$f \in L^p(0,T; L^q) \cap L^\infty (0,T ; L^\infty)$, 
$g \in L^\infty (0,T ; L^\infty)$ and $\inf  g > c $ . 

\item Let $d=2$. 
There exists $C > 0$ such that for every $f \in L^2(0,T;H^1(\mathbb T^2))$ 
\[
\|\nabla f_{\ep}\|_{L^2(0,T; L^{\infty})} \le 
C \ep^{-1} \|f\|_{L^2(0,T; H^1(\mathbb T^2))}.
\]
Moreover, for every $f \in L^2(0,T;H^1(\mathbb T^2))$ 
\begin{equation}\label{1208-3}
\limsup_{\ep \to 0} \ep \|\nabla f_{\ep}\|_{L^2(0,T; L^{\infty})} =0.
\end{equation}
\item Let $p, p_1, q, q_1 \in [1, \infty)$, $p_2 , q_2 \in(1, \infty]$, 
$1/p = 1/p_1 + p_2 , 1/q = 1/q_1 + 1/q_2$. 
Then there exists $C > 0$ such that 
for every $f \in L^{p_1}(0,T; W^{1, q_1}) $ and $ g \in L^{p_2}(0,T; L^{q_2})$ 
\[
\|(f g)_{\ep} - f_{\ep}g_{\ep}\|_{L^p(0,T; L^q)}
\le
C\ep \|f \|_{L^{p_1}(0,T; W^{1, q_1})}\|g\|_{L^{p_2}(0,T; L^{q_2})}. 
\]
Moreover,  
\begin{equation}\label{1208-4}
\limsup_{\ep \to 0} \ep^{-1} \|(fg)_{\ep} - f_{\ep}g_{\ep}\|_{L^p(0,T;L^q)} =0
\end{equation}
provided that 
$f \in L^{p_1}(0,T; W^{1, q_1}) $ and $ g \in L^{p_2}(0,T; L^{q_2})$. 
\end{enumerate}
\end{lem}

We here give a self-contained proof. 

\begin{proof}

\noindent 
(1) 
We write 
\[
\nabla f_{\ep} (t,x) = 
\ep^{-1}\int_{\mathbb T^d} f(y) (\nabla \eta)_{\ep} (x -y) dy, 
\]
and the Young inequality implies the inequality, 
\[
\| \nabla f_{\ep} \|_{L^p(0,T; L^q)} \le 
C\ep^{-1} \|f\|_{L^p(0,T; L^q)}. 
\]
Next we prove the convergence \eqref{1208-1}. 
When $f \in C^{\infty}((0,T) \times \mathbb T^d)$, 
it is possible to obtain the convergence. A density argument for the inclusion 
$C^{\infty}((0,T) \times \mathbb T^d) \subset L^p(0,T; L^q)$ 
with the inequality above 
proves the convergence for $f \in L^p(0,T; L^q)$.

\noindent (2) 
Using Leibniz rule and (1), we have 
\begin{equation}
\begin{aligned}
\label{1126-1}
\bigg\|\nabla\frac{f_{\ep}}{g_{\ep}} \bigg\|_{L^{p}(0,T;L^q)} 
&\le
\|\nabla f_{\ep}\|_{L^p(0,T; L^q)} \frac{1}{\inf |g_{\ep}|} 
+
\|f_{\ep}\|_{L^p(0,T; L^q)}
\frac{\|\nabla g_{\ep}\|_{L^{\infty}(0,T;L^{\infty})}}
{\inf |g_{\ep}|^2} 
\\
&\le 
c^{-1} \|\nabla f_{\ep} \|_{L^p(0,T;L^q)} + 
c^{-2} \|f_{\ep}\|_{L^p(0,T;L^q)} \|\nabla g_{\ep} \|_{L^{\infty}(0,T;L^{\infty})}
\\
&\le 
\ep^{-1} ( c^{-1} + c^{-2} \|g_{\ep}\|_{L^{\infty}(0,T;L^{\infty})})
\|f\|_{L^p(0,T;L^q)} 
\le
C \ep^{-1} \|f\|_{L^p(0,T;L^q)}. 
\end{aligned}
\end{equation}
To prove the convergence \eqref{1208-2}, 
we take $g_n \in C^\infty ([0,T] \times \mathbb T^d)$ such that
\[
\|g_n - g\|_{L^p(0,T;L^q)} \to 0 \quad
as \, \, \,  n \to \infty.
\]
Using the first inequality of \eqref{1126-1} and the inequality 
of (1), we get 
\[
\begin{aligned}
& \ep \bigg\|\nabla \frac{f_{\ep}}{g_{\ep}} \bigg\|
\\
&\le 
C \ep \|\nabla f_{\ep}\|_{L^p(0,T;L^q)} + 
C \ep \|(f_{\ep}\nabla) g_{\ep}\|_{L^p(0,T;L^q)}
\\
&\le 
C \ep \|\nabla f_{\ep}\|_{L^p(0,T;L^q)} + 
C \ep \|f_{\ep}\nabla ((g_n)_\ep - g_{\ep})\|_{L^p(0,T;L^q)}
+ C \ep \|f_{\ep}\nabla (g_n)_{\ep}\|_{L^p(0,T;L^q)}
\\
&\le 
C \ep \|\nabla f_{\ep}\|_{L^p(0,T;L^q)} + 
C \|f\|_{L^\infty (0,T ; L^\infty)} 
\big( \| g_n - g\|_{L^p(0,T;L^q)} + \varepsilon\|\nabla (g_n)_{\ep}\|_{L^p(0,T;L^q)}
\big) .
\end{aligned}
\]
We then fix $n$ sufficiently large, and apply the convergence of (1) 
as $\varepsilon \to 0$ to conclude the convergence \eqref{1208-2}.

\noindent 
(3) Combining assumption and the inequality of (1), we get
\[
\| \nabla f_{\ep}\|_{L^2(0,T; L^{\infty})}
\le  
\| \nabla f\|_{L^2(0,T;L^2)} 
\|\eta_{\ep}\|_{L^2} 
\le
C \ep^{-1} \|f\|_{L^2(0,T;H^1)}.
\]
The convergence \eqref{1208-3} can be proved in the same way as \eqref{1208-1}.\\
(4) We write 
\begin{equation}\label{1208-5}
\begin{split}
(f g)_{\ep} - f_{\ep} g_{\ep} =
 \big( (fg)_\ep -f_\ep g -fg_\ep + fg  \big) 
 - (f_{\ep} - f)(g_{\ep} -g),  
\end{split}
\end{equation}
\[
\begin{split}
(fg)_\ep -f_\ep g -fg_\ep + fg 
 =& 
\int_{B_{\ep}(0)} 
(f(t,x-y) - f(t, x))(g(t,x-y) - g(t,x)) \eta_{\ep}(y) dy 
\\
=: & I_\varepsilon (f,g) .
\end{split}
\]
We estimate the first term in the right-hand side of \eqref{1208-5}, 
\[
\begin{aligned}
& \| I_\varepsilon (f,g)  \|_{L^p(0,T;L^q)}
\\
&\le C \sup_{y \in B_{\ep}(0)} 
\bigg(\int^T_0 \bigg(\int_{\mathbb T^d}
\bigg((f(t,x-y) - f(t,x))(g(t,x-y) - g(t,x)) \bigg)^q dx\bigg)^{\frac{p}{q}}
dt\bigg)^{\frac{1}{p}}
\\
&\le
\ep \bigg(\int^T_0 \bigg(
\sup_{y \in B_{\ep}(0)}\bigg\|\frac{f(t,\cdot-y) - f(t,\cdot)}{|y|}\bigg\|_{L^{q_1}} 
\|g(t,\cdot-y) - g(t,\cdot)\|_{L^{q_2}}
\bigg)^p dt\bigg)^{\frac{1}{p}} 
\\
&\le 
C\ep  \|f\|_{L^{p_1}(0,T;W^{1,q_1})}
\|g\|_{L^{p_2}(0,T;L^{q_2})},
\end{aligned}
\]
and the second term, 
\[
\begin{aligned}
&\| (f_{\ep} -f) (g_{\ep} - g)\|_{L^p(0,T;L^q)}
\\
&\le
\bigg(\int^T_0\bigg(\int_{\mathbb T^d}\bigg(
\int_{B_{\ep}(0)} 
(f(t,x-y) - f(t,x)) \eta_{\ep}(y) dy 
\bigg)^{q_1}dx
\bigg)^{\frac{p_1}{q_1}}dt\bigg)^{\frac{1}{p_1}}
\|g_{\ep}-g\|_{L^{p_2}(0,T;L^{q_2})}
\\
&\le 
\ep\bigg(\int^T_0\bigg(
\sup_{y \in B_\ep (0)}\bigg\|\frac{f(t,\cdot-y) - f(t,\cdot)}{|y|}\bigg\|_{L^{q_1}}
\bigg)^{p_1}dt\bigg)^{\frac{1}{p_1}}
\|g\|_{L^{p_2}(0,T;L^{q_2})}
\\
&\le 
C\ep \|f\|_{L^{p_1}(0,T;W^{1,q_1})}
\|g\|_{L^{p_2}(0,T;L^{q_2})},
\end{aligned} 
\]
which prove the inequality of (4).

We next prove the convergence of (4). 
We take $ f_n \in C^{\infty}([0,T] \times \mathbb T^3)$ such that
\[
\|f -f_n\|_{L^{p_1}(0,T;W^{1,q_1})} \to 0 \quad 
as \quad n \to \infty.
\]
We choose $\delta >0$ arbitrarily and $\widetilde p_2, \widetilde q_2 < \infty$ 
such that 
\[
\frac{1}{p} = \frac{1}{p_1 + \delta} + \frac{1}{\widetilde p_2},
\quad 
\frac{1}{q} = \frac{1}{q_1 + \delta} + \frac{1}{\widetilde q_2} .
\]
We utilize \eqref{1208-5} and estimate $I_\ep (f,g)$ 
by the approximation of $f$ with $f_n$ and by an analogous argument to 
the above 
\[
\begin{aligned}
&\ep^{-1}\|I_\varepsilon (f,g)  \|_{L^p(0,T;L^q)}
\\
&\leq  
\ep^{-1}\|I_\varepsilon (f -f_n ,g)  \|_{L^p(0,T;L^q)}
+ \ep^{-1}\|I_\varepsilon (f_n,g)  \|_{L^p(0,T;L^q)}
\\
&\le 
C\|f - f_n\|_{L^{p_1}(0,T;W^{1,q_1})}\|g \|_{L^{p_2}(0,T;L^{q_2})}
\\
& \quad 
+ C 
\|f_n\|_{L^{p_1+ \delta}(0,T;W^{1,q_1+\delta})}
\bigg\|
 \sup_{y \in B_\ep (0)} 
 \|g(\cdot, \cdot-y) - g(\cdot,\cdot)\|_{L^{\widetilde q_2}_x}
 \bigg\|_{L^{\widetilde p_2 } (0,T)} .
\end{aligned}
\]
The first term above is independent of $\ep > 0$ and converges to zero 
as $n \to \infty$, and for fixed $n$, the second term converges to zero 
because of the continuity of the translation in 
$L^{\widetilde p_2} (0,T; L^{\widetilde q_2})$ and 
$g \in L^{\widetilde p_2} (0,T; L^{\widetilde q_2})$. 
We then obtain 
\[
\ep^{-1}\|I_\varepsilon (f,g)\|_{L^p(0,T;L^q)}
\to 0.
\]
Similarly, as for the last term $(f - f_{\ep})(g - g_{\ep})$, we have 
\[
\begin{aligned}
& \ep^{-1} \|(f - f_{\ep})(g - g_{\ep})\|_{L^p(0,T;L^q)}
\\
\leq & 
C \|f-f_n\|_{L^{p_1}(0,T;W^{1,q_1})} 
\|g\|_{L^{p_2}(0,T;L^{q_2})}
\\&+ 
C \|f_n\|_{L^{p_1+ \delta}(0,T;W^{1,q_1+\delta})} 
 \bigg\|    \sup_{y \in B_\ep (0)} \|g(\cdot, \cdot-y) - g(\cdot,\cdot)\|_{L^{\widetilde q_2}}
 \bigg\|_{L^{\widetilde p_2}(0,T)},
\end{aligned}
\]
and its convergence to zero as $\ep \to 0$.  
Therefore we conclude \eqref{1208-4}. 
\end{proof}

\section{Proof of main theorem}

For simplicity, we may consider only the case when $u \in L^4(0,T;L^4)$ 
(see Remark below Theorem~\ref{1129-2}).

We first consider the continuity equation and the momentum equation.  
For $x \in \mathbb T^d$, let test function be the mollifier $\eta_{\ep}(x - \cdot)$. 
Multiplying the test function by the continuity equation and the momentum equation, we have
\[
\begin{aligned}
\int^t_0 \bigg\{
& \bigg\< \rho (s,\cdot), \partial_s \eta_{\ep} (x - \cdot) \bigg\> +
\bigg\< (\rho u ) (s, \cdot), \nabla (\eta_{\ep}(x-\cdot))\bigg\>  \bigg\}
ds 
\\=
&\bigg\< \rho  (t,\cdot), \eta_{\ep}(x-\cdot) \bigg\> -
\bigg\< \rho_0 (\cdot), \eta_{\ep}(x-\cdot) \bigg\>,  
\\
\int^t_0 \bigg\{
&\bigg\< (\rho u) (s,\cdot), \partial_s \eta_{\ep} (x - \cdot) \bigg\> +
\bigg\< (\rho u \otimes u) (s, \cdot), \nabla (\eta_{\ep}(x-\cdot))\bigg\> 
\\
+ \bigg\< p& (\rho, \theta) (s, \cdot), {\rm div\,}(\eta_{\ep}(x-\cdot))\bigg\>
-\bigg\< \mathbb S (s,\cdot), \nabla (\eta_{\ep}(x-\cdot))\bigg\> \bigg\}
ds 
\\=
&
\bigg\< (\rho u) (t,y), \eta_{\ep}(x-\cdot) \bigg\> -
\bigg\< (\rho_0 u_0) (y), \eta_{\ep}(x-\cdot) \bigg\>,   
\end{aligned}
\]
which imply 
\[
\left\{
\begin{split}
&- \int^t_0 
{\rm div \,}(\rho u )_{\ep} (s,x) \, ds =
\rho_{\ep} (t,x) - \rho_{0,{\ep}} (x) ,
\\
& -\int^t_0 \bigg\{
{\rm div \,}(\rho u \otimes u)_{\ep} (s,x)  
+ \nabla p(\rho, \theta)_{\ep} (s,x) 
- {\rm div\,} \mathbb S_{\ep} (s,x) \bigg\} ds 
\\
& \qquad =
(\rho u )_{\ep} (t,x) - (\rho_0 u_0)_{\ep} (x).
\end{split}
\right.
\]
Differentiating by $t$ gives 
\begin{equation}
\left\{
\begin{split}
\label{1127-4}
&\partial_t \rho_{\ep} + 
{\rm div\,}(\rho u)_{\ep} =0,
\\
&\partial_t (\rho u)_{\ep} + 
{\rm div\,}(\rho u \otimes u)_{\ep} + 
\nabla p(\rho, \theta)_{\ep} - 
{\rm div \,} \mathbb S_{\ep} =0, 
\end{split}
\right.
\end{equation}
for all $(t,x) \in (0,T) \times \mathbb T^d$. 
The multiplication by $\rho^{-1}_{\ep}(\rho u)_{\ep}$ 
of the second equation of \eqref{1127-4} yields that 
\begin{equation}\label{1214-2}
\begin{aligned}
\displaystyle
0= &\int^t_s \bigg\{ 
\bigg\< \partial_{\tau} (\rho u)_{\ep},
\frac{(\rho u)_{\ep}}{\rho_{\ep}} \bigg\> 
+
\bigg\< \div (\rho u \otimes u)_{\ep}, 
\frac{(\rho u)_{\ep}}{\rho_{\ep}} \bigg\> 
\\
& \qquad +
\bigg\< \nabla p(\rho, \theta)_{\ep}, \frac{(\rho u)_{\ep}}{\rho_{\ep}} \bigg\> - 
\bigg\< {\rm div \,} \mathbb S_\varepsilon, 
\frac{(\rho u)_{\ep}}{\rho_{\ep}} \bigg\>  
\bigg\} d\tau
\\
=:
& (A)+(B)+(C)+(D) .
\end{aligned}
\end{equation}
We here extract the important terms by the lemma below.  
\begin{lem}\label{lem:0105}
The equality \eqref{1214-2} is equivalent to 
\begin{equation}
\begin{aligned}\label{1214-1}
& \frac{1}{2}\int^t_s \int_{\mathbb T^d}
\partial_{\tau} \bigg\{ \frac{(\rho u)^2_{\ep}}{\rho_{\ep}}  \bigg\} dx d\tau +
\int^t_s \int_{\mathbb T^d}
\partial_{\tau} ( \rho_{\ep} h(\rho_{\ep})) dx d\tau 
\\
&- \int^t_s \< \{ p_{th} (\rho, \theta)\}_{\ep}, {\rm div \,} u_{\ep} \> d\tau + \int^t_s \< \mathbb S_{\ep} : \nabla u_{\ep} \> d\tau 
+ R_\varepsilon (t,s)=0,
\end{aligned}
\end{equation}
where $h$ is defined by \eqref{0202-2}, 
and the error term $R_\varepsilon (t,s)$ satisfies 
\[
\lim_{\varepsilon \to 0}\sup_{s \in (0,t)} | R_\ep (t,s)| = 0 \text{ for all } t. 
\]
\end{lem}

\begin{pf} 
We write that 
\begin{equation}\label{1214-4}
\begin{split}
(A) = 
&\int^t_s
\bigg\< \partial_{\tau} (\rho u)_{\ep}, \frac{(\rho u)_{\ep}}{\rho_{\ep}} \bigg\> d\tau
\\
= 
&\frac{1}{2} \int^t_s \int_{\mathbb T^d}
 \partial_{\tau} \bigg\{ \frac{(\rho u)^2_{\ep}}{\rho_{\ep}}\bigg\}  dxd\tau + 
\frac{1}{2} \int^t_s
\bigg\< \partial_{\tau}\rho _{\ep},  \frac{(\rho u)^2_{\ep}}{\rho^2_{\ep}} \bigg\> d\tau.
\end{split}
\end{equation}
By integration by parts, 
\[
\begin{aligned}
(B) =&
-\int^t_s
\bigg\< (\rho u \otimes u)_{\ep}, \nabla \frac{(\rho u)_{\ep}}{\rho_{\ep}} \bigg\> d\tau
\\
=&  
-\int^t_s
\bigg\< (\rho u \otimes u)_{\ep} -  u_{\ep} \otimes (\rho u)_{\ep}, \nabla \frac{(\rho u)_{\ep}}{\rho_{\ep}} \bigg\>d\tau
-\int^t_s
\bigg\< u_{\ep} \otimes (\rho u)_{\ep}, \nabla \frac{(\rho u)_{\ep}}{\rho_{\ep}} \bigg\> d\tau 
\\
=: &
(B_1) + (B_2),
\end{aligned}
\]
and 
\[
\begin{aligned}
(B_2) 
=&\int^t_s \bigg\< {\rm div\,} \Big( u_{\ep} \otimes (\rho u)_{\ep}\Big) , \frac{(\rho u)_{\ep}}{\rho_{\ep}} \bigg\> d\tau 
\\
=&
 \int^t_s \bigg\< {\rm div}\, u_{\ep}, \frac{(\rho u)^2_{\ep}}{\rho_{\ep}} \bigg\> d\tau
+\frac{1}{2}\int^t_s \bigg\< \frac{u_{\ep}}{\rho_{\ep}}, \nabla 
(\rho u)^2_{\ep} \bigg\> d\tau
\\
=&
\frac{1}{2}
\int^t_s \bigg\{\bigg\< \rho_{\ep}{\rm div\,}u_{\ep}, \frac{(\rho u)^2_{\ep}}{\rho^2_{\ep}} \bigg\> 
- \bigg\< u_{\ep} \nabla \rho_{\ep}, \frac{(\rho u)^2_{\ep}} {\rho^2_{\ep}} \bigg\> \bigg\}d\tau
\\
=&
 \frac{1}{2}\int^t_s \bigg\< {\rm div}\,(\rho_{\ep} u_{\ep}), \frac{(\rho u)^2_{\ep}} {\rho^2_{\ep}} \bigg\> d\tau
\\
=&
\frac{1}{2}\int^t_s \bigg\< {\rm div}\,\{(\rho_{\ep} u_{\ep})- (\rho u)_{\ep} \}, \frac{(\rho u)^2_{\ep}} {\rho^2_{\ep}} \bigg\> d\tau
+
\frac{1}{2}\int^t_s \bigg\< {\rm div}\,(\rho u)_{\ep} , \frac{(\rho u)^2_{\ep}} {\rho^2_{\ep}} \bigg\> d\tau
\\
=&
 (B_{2,1}) - 
\frac{1}{2}\int^t_s \bigg\< \partial_{\tau} \rho_{\ep}, \frac{(\rho u)^2_{\ep}} {\rho^2_{\ep}} \bigg\> d\tau.
\end{aligned}
\]
We notice that the second term above line is cancelled with 
the last term in \eqref{1214-4}. 
We prove that $(B_1), (B_{2,1})$ tends to $0$ as $\varepsilon \to 0$. 
By H\"older's inequality, 
\[
\begin{aligned}
|(B_1)| = &\bigg|
\int^t_s \bigg\< (\rho u \otimes u)_{\ep} - u_{\ep} \otimes (\rho u)_{\ep}, \nabla \frac{(\rho u)_{\ep}}{\rho_{\ep}} \bigg\> d\tau \bigg|
\\
\le 
&\bigg\| \nabla \frac{(\rho u)_{\ep}}{\rho_{\ep}}\bigg\|_{L^4(0,T ;L^4)}
\|(\rho u \otimes u)_{\ep} - u_{\ep} \otimes (\rho u)_{\ep}\|_{L^{\frac{4}{3}}(0,T; L^{\frac{4}{3}})}. 
\end{aligned}
\]
By the assumption $u \in L^4(0,T;L^4)$, 
we can applying the inequality of (2) and the convergence of (4) in Lemma~\ref{0909-1}, 
and obtain 
\[
\begin{split}
&\bigg\|\nabla \frac{(\rho u)_{\ep}}{\rho_{\ep}} \bigg\|_{L^4(0,T; L^4)} \le
 C\ep^{-1} \|u\|_{L^4(0,T; L^4)},
 \\
&\limsup_{\ep \to 0}\ep^{-1} \|(\rho u \otimes u)_{\ep} - u_{\ep} \otimes (\rho u)_{\ep}\|_{L^{\frac{4}{3}}(0,T; L^{\frac{4}{3}})} =0,
\end{split}
\]
which yield 
\[
\limsup_{\ep \to 0} \sup_{s \in (0,t)}|(B_1)| = 0.
\]
As for $(B_{2,1})$, we use integration by parts, then we have
\[
\begin{aligned}
|(B_{2,1})| &= \frac{1}{2} \bigg| \int^t_s 
\bigg\< (\rho u )_{\ep} - \rho_{\ep}u_{\ep}, \nabla \frac{(\rho u)^2_{\ep}}{\rho^2_{\ep}} \bigg\> d\tau
\bigg|  \\
&\le 
C \|(\rho u )_{\ep} - \rho_{\ep}u_{\ep}\|_{L^2(0,T; L^2)} 
\bigg\|\nabla \frac{(\rho u)^2_{\ep}}{\rho^2_{\ep}}\bigg\|_{L^2(0,T; L^2)} .
\end{aligned}
\]
We again apply the inequality of (2) in Lemma~\ref{0909-1} 
 and the convergence of (4) in Lemma~\ref{0909-1}, 
 and then obtain 
\[
\bigg\| {\rm div}\, \frac{|(\rho u)_{\ep}|^2}{\rho_{\ep}^2}\bigg\|_{L^2(0,T;L^2)} \le 
C \ep^{-1} \| u \|^2_{L^4(0,T; L^4)}, 
\]
\[
\limsup_{\ep \to 0} \ep^{-1} \|(\rho u)_{\ep} - \rho_{\ep}u_{\ep} \|_{L^2(0,T;L^2)} =0,  
\]
since $\rho \in  L^{\infty}(0,T;L^{\infty}) $, 
$u \in L^2 (0,T ; W^{1,2})$.
Therefore we extract the important part of $(A) + (B)$ with 
the error estimate that 
\[
\lim_{\varepsilon\to0} 
\sup _{s \in (0,t)}\left| (A) + (B) - \frac{1}{2} \int^t_s \int_{\mathbb T^d}
 \partial_{\tau} \bigg\{ \frac{(\rho u)^2_{\ep}}{\rho_{\ep}}\bigg\}  dxd\tau
 \right| = 0 .
\]

We turn to consider $(C)$ and write 
\[
(C) = \int^t_s \bigg\{
\bigg\< \nabla (p_e(\rho))_{\ep} ,  \frac{(\rho u)_{\ep}}{\rho_{\ep}} \bigg\>
+ \bigg\< \nabla (p_{th}(\rho, \theta))_{\ep},  \frac{(\rho u)_{\ep}}{\rho_{\ep}} \bigg\> \bigg\} d\tau
 = (C_1) + (C_2). 
\]
We extract the second term and the third term of \eqref{1214-1} 
from $(C_1), (C_2)$, respectively. 
We first estimate $(C_1)$ and write 
\[
\begin{aligned}
(C_1) &=  \int^t_s
\bigg\< \nabla (p_e(\rho))_{\ep} , \frac{(\rho u)_{\ep}}{\rho_{\ep}} \bigg\> 
d\tau
\\
&=
\int^t_s\bigg\{
\bigg\< \nabla \bigg\{(p_e(\rho))_{\ep} -  p_e(\rho_{\ep})\bigg\}, \frac{(\rho u)_{\ep}}{\rho_{\ep}} \bigg\>
+ \bigg\< \frac{\nabla p_e(\rho_{\ep})}{\rho_{\ep}}, (\rho u)_{\ep} \bigg\>
\bigg\} d\tau 
\\
&= (C_{1,1}) + (C_{1,2}).
\end{aligned}
\]
We can see that $(C_{1,2})$ becomes the second term of \eqref{1214-1}.  
In fact, we introduce 
\[
\displaystyle
g (\rho_\varepsilon) := \int^{\rho_\varepsilon}_1 \frac{p'_e(z)}{z} dz ,
\]
and notice that 
\[
\frac{\nabla p_e(\rho_{\ep})}{\rho_{\ep}} 
= \nabla g(\rho_\varepsilon), 
\]
and by integration by parts that 
\[
\begin{split} 
g (\rho_\varepsilon) 
=& h(\rho_\varepsilon) + \frac{p_e(\rho_\varepsilon)}{\rho_\varepsilon} - p_e (1), 
\\
( \partial _t \rho_{\varepsilon} ) g(\rho_\varepsilon) 
=& \partial _t ( \rho_{\varepsilon} g(\rho_\varepsilon))
  -  \rho_{\varepsilon}\partial _t g(\rho_\varepsilon)
\\
=&  \partial _t 
 \Big(\rho_{\varepsilon} \Big( h(\rho_\varepsilon) + \frac{p_e(\rho_\varepsilon)}{\rho_\varepsilon} - p_e (1)
 \Big)\Big)  
- \rho_\varepsilon \frac{p_e'(\rho_\varepsilon) }{\rho_\varepsilon}
  \partial_t \rho_\varepsilon
\\
=& \partial _t (\rho_\varepsilon h(\rho_\varepsilon)) 
- p_e (1) \partial _t \rho_\varepsilon 
\\
=& \partial _t (\rho_\varepsilon h(\rho_\varepsilon)) 
+ p_e (1) {\rm div \, } ( \rho u)_\varepsilon ,
\end{split}
\]
which imply that 
\[
\begin{aligned}
(C_{1,2})
&= \int_s^t \bigg\< \nabla g(\rho_\ep), (\rho u)_{\ep} \bigg\>
= - \int_s^t \bigg\< g(\rho_\ep), {\rm div \, }(\rho u)_{\ep} \bigg\>
=  \int_s^t \bigg\< g(\rho_\ep), \partial _\tau \rho_{\ep} \bigg\>
d\tau 
\\
&= \int_s^t \int_{\mathbb T^d} 
 \Big( \partial_\tau ( \rho_{\ep} h(\rho_{\ep}) ) + 
 p_e(1) {\rm div}\, (\rho u)_{\ep} \Big) 
 ~dxd\tau
\\
& = \int_s^t \int_{\mathbb T^d}  \partial_\tau ( \rho_{\ep} h(\rho_{\ep}) ) 
~dxd\tau. 
\end{aligned}
\]
We next show that $(C_{1,1})$ converges to $0$. 
We approximate $(\rho u)_\varepsilon$ by $\rho_\varepsilon u_\varepsilon$ 
and apply integration by parts, and then have that 
\begin{equation}
\begin{aligned}
\label{0915-1}
|(C_{1,1})|
\leq &
\bigg|\int^t_s 
\bigg\< \nabla \bigg\{(p_e(\rho))_{\ep} - p_e(\rho_{\ep})\bigg\}, \frac{(\rho u)_{\ep} - \rho_{\ep}u_{\ep}}{\rho_{\ep}} 
\bigg\>
d\tau \bigg|
\\
& + 
\bigg|
\int^t_s \bigg\< \bigg\{(p_e(\rho))_{\ep} - p_e(\rho_{\ep})\bigg\}, {\rm div} \, u_{\ep} \bigg\> d\tau
\bigg| \\
\leq &
\bigg\| \nabla \bigg\{(p_e(\rho))_{\ep} - p_e(\rho_{\ep})\bigg\} \bigg\|_{L^{\infty}(0,T ;L^{\infty})}
\bigg\| \frac{(\rho u)_{\ep} - \rho_{\ep}u_{\ep}}{\rho_{\ep}}\bigg\|_{L^1(0,T; L^1)} 
\\ 
 & + 
\| (p_e(\rho))_{\ep} - p_e(\rho_{\ep}) \|_{\bl L^2(0,T; L^2)}
\|{\rm div \,} u_{\ep} \|_{\bl L^2(0,T; L^2)}.
\end{aligned}
\end{equation}


We  deal with the first term of right-hand side of \eqref{0915-1}. 
By Lemma \ref{0909-1} (1), we get
\[
\begin{aligned}
\bigg\| \nabla \bigg\{(p_e(\rho))_{\ep} - p_e(\rho_{\ep})\bigg\} \bigg\|_{L^{\infty}(0,T; L^{\infty})} &\le
\| \nabla (p_e(\rho))_{\ep} \|_{L^{\infty}(0,T;L^{\infty})} + \| p_e'(\rho) \nabla \rho_{\ep} \|_{L^{\infty}(0,T;L^{\infty})} 
\\
&\leq 
 \ep^{-1}(\|p_e(\rho)\|_{L^{\infty}(0,T; L^{\infty})} +
 \| \rho \|_{L^{\infty}(0,T;L^{\infty})})
 \le C \ep^{-1}
\end{aligned}
\]
with $ \rho \in L^{\infty}(0,T; L^{\infty}(\mathbb T^2))$ and the assumption \eqref{0915-3}.
Since $u \in L^2(0,T; W^{1,2})$ and $ \rho \in L^{\infty}(0,T; L^{\infty}) \subset L^2(0,T;L^2)$, 
we get by Lemma \ref{0909-1} (4), 
\[
\limsup_{\ep \to 0} \ep^{-1} \bigg\|\frac{(\rho u)_{\ep} - \rho_{\ep}u_{\ep}}{\rho_{\ep}} \bigg\|_{L^1(0,T ; L^1)} \le
C \limsup_{\ep \to 0} \ep^{-1} \|(\rho u)_{\ep} - \rho_{\ep}u_{\ep} \|_{L^1(0,T; L^1)} = 0.
\]
It follows that 
\begin{equation}
\label{1223-1}
\limsup_{\ep \to 0} 
\bigg\| \nabla \bigg\{(p_e(\rho))_{\ep} - p_e(\rho_{\ep})\bigg\} \bigg\|_{L^{\infty}(0,T;L^{\infty})} \bigg\| \frac{(\rho u)_{\ep} - \rho_{\ep}u_{\ep}}{\rho_{\ep}} \bigg\|_{L^1(0,T; L^1)} =0. 
\end{equation}
{\bl For the second term of the right-hand side of \eqref{0915-1} 
it suffices to show that
\begin{equation}
\label{0915-2}
\|(p(\rho))_{\ep} - p(\rho_{\ep})\|_{L^2(0,T; L^2)} \to 0.
\end{equation}
By $c_1 \leq \rho \leq c_2$ from the assumption \eqref{0915-3},  
the mean value theorem 
and $\eta _\ep$ having the unit mass, we have 
\begin{equation}
\begin{aligned} \notag 
&| (p_e(\rho))_{\ep}(\tau, x) - p_e(\rho_{\ep}) (\tau, x)| \\
&\le
\bigg|
 \int_{\mathbb T^d} \{ p_e (\rho(\tau, x-y))
\}\eta_{\ep} (y) dy
- p_e(\rho(\tau, x))
\bigg|
\\
&\quad+
\bigg| p_e(\rho(\tau, x)) - p_e\bigg( \int_{\mathbb T^d}\rho(\tau, x-y) \eta_{\ep}(y) dy \bigg) 
 \bigg| 
\\
&\le 
2 \sup_{c \in (c_1, c_2)} |p'_e(c)|
\bigg|
\int_{\mathbb T^d} 
 \Big( \rho(\tau, x-y) - \rho(\tau, x)  \Big) 
 \eta_{\ep} (y)
 dy \bigg|
.
\end{aligned}
\end{equation}
By $\rho \in L^{\infty}(0,T; L^{\infty}) \subset L^2(0,T;L^2)$ and the continuity of the translation in $L^2(0,T;L^2)$, we have
\[
\begin{aligned}
&\limsup_{\ep \to 0}
\| (p_e(\rho))_{\ep} - p_e(\rho_{\ep}) \|_{L^2(0,T;L^2)}
\\
\le
&2\limsup_{\ep \to 0}
 \sup_{c \in [c_1, c_2]} |p'_e(c)|
\bigg(\int^T_0 
\sup_{y \in B_{\ep}(0)} 
\|\rho(\tau, \cdot-y) - \rho(\tau, \cdot) \|^2_{L^2(\mathbb T^d)} d\tau \bigg)^{\frac{1}{2}} =0,
\end{aligned}
\]
which proves \eqref{0915-2}. 
}
We then conclude by \eqref{1223-1} and \eqref{0915-2} that 
\[
\limsup_{\ep \to 0} \sup_{s \in (0,t)} |(C_{1,1})| =0 . 
\]
As for $(C_2)$, we write 
\[
\begin{aligned}
(C_2) 
&= 
\int^t_s \bigg\{ \bigg\< \nabla\{p_{th}(\rho, \theta)\}_{\ep}, \frac{(\rho u)_{\ep} - \rho_{\ep} u_{\ep}}{\rho_{\ep}} \bigg\> 
+ \bigg\< \nabla\{p_{th}(\rho, \theta)\}_{\ep},  u_{\ep} \bigg\> \bigg\} d\tau .
\end{aligned}
\]
The second term of the right hand side is nothing but the third term of 
\eqref{1214-1}, and we show the first term above coverges to $0$ as 
$\ep \to 0$. 
By H\"older's inequality, Lemma~\ref{0909-1} (1), and the assumption of $\rho \in L^{\infty}(0,T; L^{\infty})$, we get
\[
\begin{aligned}
& 
\left| 
\int^t_s 
\bigg\< \nabla\{p_{th}(\rho, \theta)\}_{\ep}, \frac{(\rho u)_{\ep} - \rho_{\ep} u_{\ep}}{\rho_{\ep}} \bigg\> 
d\tau 
\right| 
\\
 &\le 
\| \nabla \{p_{th}(\rho, \theta)\}_{\ep} \|_{L^2(0,T;L^2)}
\bigg\| \frac{(\rho u)_{\ep} - \rho_{\ep}u_{\ep}}{\rho_{\ep}} \bigg\|_{L^2(0,T; L^2)}
\\ & \le C \ep^{-1} \|  p_{th}(\rho, \theta) \|_{L^2(0,T; L^2)} 
\|(\rho u)_{\ep} - \rho_{\ep}u_{\ep}\|_{L^2(0,T;L^2)}. 
\end{aligned}
\]
It follows from
 $\rho \in L^{\infty}(0,T; L^{\infty})$, $u \in L^2(0,T; W^{1,2})$ 
 and Lemma~\ref{0909-1} (4) that 
\[
\limsup_{\ep \to 0} \ep^{-1} 
\|(\rho u)_{\ep} - \rho_{\ep}u_{\ep}\|_{L^2(0,T;L^2)} = 0.
\]
Therefore, we conclude 
\[
\limsup_{\ep \to 0} \sup_{s \in (0,t)}
\left| (C) 
- \int_s^t \int_{\mathbb T^d}  
\Big( \partial_t ( \rho_{\ep} h(\rho_{\ep}) ) 
- \{ p_{th} (\rho, \theta)\} _\ep  \, {\rm div \, } u_\ep 
\Big)
~dxd\tau
\right| 
 =0.
\]

Finally, we consider $(D)$. 
We write 
\[
\begin{aligned}
(D) &= 
 -\int^t_s \bigg\< {\rm div} \, \mathbb S_{\ep}, \frac{(\rho u)_{\ep}}{\rho_{\ep}} \bigg\> d\tau 
= 
-\int^t_s \bigg\{ \bigg\< {\rm div} \, \mathbb S_{\ep}, \frac{(\rho u)_{\ep} - \rho_{\ep} u_{\ep}}{\rho_{\ep}} \bigg\> 
+ \bigg\< {\rm div } \, \mathbb S_{\ep},  u_{\ep} \bigg\> \bigg\} d\tau . 
\end{aligned}
\]
By H\"older's inequality, Lemma~\ref{0909-1} (1) and assumption $\rho \in L^{\infty}(0,T;L^{\infty})$, 
\[
\begin{aligned}
\left| 
\int_s^t \bigg\< {\rm div} \, \mathbb S_{\ep}, \frac{(\rho u)_{\ep} 
  - \rho_{\ep} u_{\ep}}{\rho_{\ep}} \bigg\>
  d\tau 
\right|  &\le 
\| {\rm div} \, \mathbb S_{\ep} \|_{L^2(0,T ;L^2)}
\bigg\| \frac{(\rho u)_{\ep} - \rho_{\ep}u_{\ep}}{\rho_{\ep}} \bigg\|_{L^2(0,T;L^2)}
\\ & \le C \ep^{-1} \| \mathbb S \|_{L^2(0,T;L^2)} 
\|(\rho u)_{\ep} - \rho_{\ep}u_{\ep}\|_{L^2(0,T;L^2)}.
\end{aligned}
\]
Since $\rho \in L^{\infty}(0,T;L^{\infty})$ and $u \in L^2(0,T; W^{1,2})$, by Lemma~\ref{0909-1} (4), we have 
\[
\limsup_{\ep \to 0} \ep^{-1} 
\|(\rho u)_{\ep} - \rho_{\ep}u_{\ep}\|_{L^2(0,T;L^2)} = 0,
\]
and the convergence of the first term. 
We then conclude
\[
\limsup_{\ep \to 0} \sup_{s \in (0,t)}
\left|(D) -
\int_s^t 
 \bigg\< \mathbb S_{\ep}:  \nabla u_{\ep} \bigg\> \bigg\} d\tau 
 \right| =0.
\]

Therefore we complete the proof of \eqref{1214-1} 
with the error estimate. 
\end{pf}

\vskip3mm

We prove Theorem~\ref{1129-2} with the help of Lemma~\ref{lem:0105}.  
We consider the thermal	equation by chosing a test function 
as a constant function. It follows from the definition of the 
weak solution and $\partial _t 1 = \partial _x 1 = 0$ that 
\begin{equation}\label{1012-2}
\begin{aligned}
&\int_{\mathbb T^d}\rho Q(\theta)(t,x) dx
 =
&\int_{\mathbb T^d} \rho Q(\theta)(s,x) dx + 
\int^t_s \bigg\{\< \mathbb S : \nabla u \> - 
\< p_{th}(\rho, \theta)\, ,{\rm div\,} u \> \bigg\}d\tau.
\end{aligned}
\end{equation}
By adding \eqref{1214-1} to \eqref{1012-2}, we obtain
\begin{equation}\label{0105-2}
\begin{aligned}
&
\frac{1}{2}\int_{\mathbb T^d}
\bigg\{\frac{(\rho u)^2_{\ep}}{\rho_{\ep}} \bigg\}(t,x)dx + 
\int_{\mathbb T^d}\rho_{\ep} h (\rho_{\ep}) (t,x) dx + 
\int_{\mathbb T^d}\rho Q(\theta)(t,x)dx
\\
=
&\frac{1}{2}\int_{\mathbb T^d}
\bigg\{\frac{(\rho u)^2_{\ep}}{\rho_{\ep}} \bigg\}(s,x) dx + 
\int_{\mathbb T^d}\rho_{\ep} h (\rho_{\ep}) (s,x) dx + 
\int_{\mathbb T^d}\rho Q(\theta)(s,x)dx   
\\&
+ \int^t_s 
\< \{p_{th}(\rho, \th)\}_{\ep},{\rm div \,}u_{\ep} \> d\tau
- \int^t_s \< p_{th}(\rho, \theta),  {\rm div \,} u\> d\tau 
\\
&- \int^t_s \<\, \mathbb S_{\ep} : \nabla u_{\ep} \>  d\tau 
+ \int^t_s \< \mathbb S : \nabla u \> d\tau 
- R_{\ep}(t,s).
\end{aligned}
\end{equation}
We start by taking the limit as $s \to 0$ for each $\ep > 0$. 
The weak continuity of $\rho u$ gives the pointwise convergence of 
$(\rho u)_\ep$ to $\rho u$, and the Lebesgue dominated convergence theorem 
implies that 
\[
\frac{1}{2} \int_{\mathbb T^d}
\frac{(\rho u)^2_{\ep}}{\rho_{\ep}}(s,x) dx \to 
\frac{1}{2} \int_{\mathbb T^d}
\frac{(\rho_0 u_0)^2_{\ep}}{(\rho_0)_{\ep}}(x) dx , 
\quad \text{as } s \to 0 \text{ for each } \ep > 0 .
\]
Similarly, we also have from the weak continuity of $\rho$ that 
\[
\int_{\mathbb T^d} \rho_{\ep}h(\rho_{\ep}) (s,x)dx 
\to 
\int_{\mathbb T^d} (\rho_0)_{\ep}h((\rho_0)_{\ep}) (x)dx ,
\quad \text{as } s \to 0 \text{ for each } \ep > 0 .
\]
The convergence of the third term $\rho Q(\theta)(s,x)$ as $s \to 0$ 
follows from the weak continuity of $\rho Q(\theta)$ due to the definition 
of the weak solutions. As for the integrals, the well-definedness is 
assured by the definition of the weak solution 
and it is possible to take the limit as $s \to 0$ due to the integrability, 
and we will apply Lemma~\ref{lem:0105} to the error term $R_\ep(t,s)$. 
Finally, we take the limit as $\ep \to 0$. The integrability 
of $\rho u^2, \rho h (\rho), p_{th}(\rho, \theta),  {\rm div \, }u, 
\mathbb S, \nabla u$ and an elemental property of the mollifier imply that 
\[
\begin{aligned}
& \frac{1}{2}\int_{\mathbb T^d}
\frac{(\rho u)^2_{\ep}}{\rho_{\ep}}(\tau,x) dx
\to 
\frac{1}{2} \int_{\mathbb T^d}
\rho u ^2(\tau,x) dx , 
\quad \text{ for } \tau = 0 ,t,
\\
& \int_{\mathbb T^d}\rho_{\ep} h(\rho_{\ep}) (\tau,x)dx
\to
\int_{\mathbb T^d}\rho h(\rho)(\tau,x) dx,
\quad \text{ for } \tau = 0 ,t, 
\\
& \int^t_0 \< \{p_{th}(\rho, \theta)\}_{\ep}, {\rm div \,}u_{\ep}\> d\tau - 
\int^t_0 \< (p_{th}(\rho, \theta),  {\rm div \,}u)\> d\tau \to 0, 
\\
& \int^t_0 \<\, \mathbb S_{\ep} : \nabla u_{\ep} \>  d\tau 
- \int^t_0 \< \mathbb S : \nabla u \> d\tau \to 0,
\end{aligned}
\]
as $\ep \to 0$, which proves the energy equality in Theorem~\ref{1129-2}.

\end{document}